\documentclass[12pt]{article}
\usepackage{amsmath,amsthm}
\usepackage[a4paper]{geometry}

\theoremstyle{definition}
\theoremstyle{plain}
\newtheorem{defn}{Definition}
\newtheorem{thm}{Theorem}
\newcommand{\real}{\mbox{$I\!\!R$}}

\date{}

\title{Exterior boundary-value Poincaré problem for elliptic systems
  of the second order with two independent variables}

\author{ F. Criado-Aldeanueva\footnote{Department of Applied Physics
  II, Polytechnic School, Malaga University. E-mail:
  fcaldeanueva@ctima.uma.es}
  \and N. Odishelidze\footnote{Department of Computer Sciences,
  Faculty of Exact and Natural Sciences, Iv. Javakhishvili Tbilisi
  State University. E-mail: nana\_georgiana@yahoo.com}
  \and
  J. M. Sanchez\footnote{Department of didactic of Mathematics, Malaga
    University. E-mail: jmss@uma.es}
  \and M. Khachidze\footnote{Computer Sciences Department,
    Iv. Javakhishvili Tbilisi State University, Tbilisi, Georgia.
    E-mail: manana.khachidze@tsu.ge}
}

\begin{document}

\maketitle

\begin{abstract}
  This paper offers a number of examples showing that in the case of
  two independent variables the uniform ellipticity of a linear system
  of differential equations with partial derivatives of the second
  order, which fulfills condition (\ref{eq:3}), do not always cause
  the normal solvability of formulated exterior elliptic problems in
  the sense of Noether.

  Nevertheless, from the system of differential equations with partial
  derivatives of elliptic type it is possible to choose, under certain
  additional conditions, classes which are normally solvable in the
  sense of Noether.

  This paper also shows that for the so-called decomposed system of
  differential equations, with partial derivatives of an elliptic type
  in the case of exterior regions, the Noether theorems are valid.

    \vskip
  0.5cm
  \noindent {\bf MSC Classification (2010):} 35J57,
  49K20, 49K21.
\end{abstract}

\section{Introduction}

Among the linear elliptic boundary-value problems, Poincare problem is
of great importance, side by side Dirichlet and Neumann problems.
Numerous processes occurring in the continuum (for example, sea tides,
see \cite{13}, as well as \cite{12}) can be simulated in terms of this
problem. This problem differs essentially from Dirichlet and Neumann
problem in the fact that it is normally solvable according to Noether
under rather general assumptions (i.e. Noether theorems known in the
theory of singular integral equations are valid for it, see
\cite{13,15,16,21,17,18,20,14,19}).

In the process of the investigation of Dirichlet and Neumann problems
for external regions, certain difficulties arose due to the complex,
behavior of the solutions of elliptic equations at infinity. This fact
was revealed on the example of Helmholtz equation simulating wave
processes in the linear formulation. The irradiation principle is of
great importance for the theory of wave processes investigated in
physics, technology, ecology and natural science. The existence,
uniqueness and stability of the solutions of the mentioned problems
for Helmholtz equation in infinite regions with boundary components
representing closed Lyapunov's surfaces (in three-dimensional case)
and closed Lyapunov's curves (in two-dimensional case) were
established in the class of functions satisfying the condition of
Sommerfeld's irradiation at infinity \cite{22,23,24}.

\section{Statements and definitions}

Let $S$ be a closed Lyapunov curve on the plane of complex variable
$z= x + i y$, and $D^-$ a simply connected domain with boundary $S$,
containing an infinitely distant point of the plane (exterior domain).

This paper is devoted to the study of the exterior boundary-value
Poincaré problem for an elliptic system of the type
\begin{equation}
  A \frac{\partial^2 u}{\partial x^2} + 2 B \frac{\partial^2 u}{\partial x \partial y} + C \frac{\partial^2 u}{\partial y^2} = 0, \label{eq:1}
\end{equation}
where $A = (a_{jk})$, $B = (b_{jk})$, $C = (c_{jk})$ are given real $n
\times n$ constant matrices and $w = (u_1, \ldots, u_n)$ the search
real vector.

\begin{defn}\label{def:1}
  Vector $w(x,y)$, having continuous derivatives of the second order
  bounded at infinity and satisfying equation (\ref{eq:1}) in domain
  $D^-$, will be called the regular solution of this equation.
\end{defn}

Let real $n \times n$ matrices $P=(p_{jk})$, $Q=(q_{jk})$,
$R=(r_{jk})$ and the real vector $f = (f_1,\ldots,f_n)$ be given on
$S$. Since $S$ is a rectifiable curve, it is assumed that $P$, $Q$,
$R$ and $f$, are given as the functions of arc length $s$ of the curve
reckoned from the fixed point on $S$ in the direction that domain
$D^-$ leaves on the right.

\begin{defn}\label{def:2}
  The exterior boundary-value Poincaré problem is assumed as the
  problem of determination of regular solution $w(z) = w(x,y)$ in
  domain $D^-$ of equation (\ref{eq:1}), satisfying the boundary-value
  condition on $S$
  \begin{equation}
    P \frac{\partial u}{\partial x} + Q \frac{\partial u}{\partial y}
    + R w = f. \label{eq:2}
  \end{equation}
\end{defn}

It is assumed that the finite value of vector $w(z)$ and its first
derivatives inside $D^-$ on $S$ exist. Note that when $n=1$ and
$R\equiv 0$, problem (\ref{eq:1}), (\ref{eq:2}) is called an oblique
derivative problem. Since in each one of the $n$ boundary-value
conditions (\ref{eq:2}) there are oblique derivatives along different
directions for various components of search bounded vector $w =
(u_1,\ldots,u_n)$ in $D^-$, it was necessary to call problem
(\ref{eq:2}) that of Poincaré. At $p_{jk} = q_{jk} = 0$, $j \neq k$,
$p_{jj} = \cos \widehat{\nu x}$, $q_{jj} = \cos \widehat{\nu y}$, $j,k
= 1, \ldots,n$, where $\nu$ inner normal to $S$ at the point $(x,y)$
relative to $D^-$, boundary-value condition (\ref{eq:2}) is somewhat
different to that of Neumann. If $P \equiv Q \equiv 0$ and $\det(R)
\neq 0$ on $S$, condition (\ref{eq:2}) is transferred into the
boundary-value condition of the Dirichlet problem.

The case in which the roots of the characteristic determinant of
system (\ref{eq:1}) are complex will be considered throughout this
paper, such system is called elliptic \cite{2}. At $n=1$ equality
(\ref{eq:1}) represents an elliptic equation of the second order,
which is reduced to Laplace's equation by means of non-singular
(non-degenerate) affine transformation of $x,y$ variables.


\section{Examples}

The domain $D^-$ is defined as the exterior of the unit circle $|z|
\le 1$ for the examples only in this section.


For Laplace's equation, assuming that $S$ is Lyapunov curve, $P$, $Q$,
$R$, $f$ are H\"older continuous \cite{11} on the whole $S$ i.e. $p,
q, r, f \in C^{0,h}(S)$, and
\begin{equation}
  \det (P + i Q) \neq 0, \label{eq:3}
\end{equation} 
for problem (\ref{eq:2}) the normal solvability according to Noether
holds \cite{1,2}, as well as the uniqueness of the solution to
Dirichlet problem.

In the case of the elliptic system of type (\ref{eq:1}), the
investigation of exterior problems of Dirichlet, Neumann and Poincaré
becomes rather difficult. Actually, at $n > 1$ as in the case of
bounded domains \cite{1}, the requirement of uniform ellipticity of
system (\ref{eq:1}) does not always guarantee the normal solvability
(accordingly to Fredholm and Noether) of exterior problem
(\ref{eq:2}). This is easy to check by the example of Bidsatze
elliptic system\footnote{
  System (\ref{eq:4}) can be written as:
  \begin{equation*}
    \begin{pmatrix}1 & 0\\0 & 1\end{pmatrix}
      \begin{pmatrix}\frac{\partial^2 u_1}{\partial x^2}\\\frac{\partial^2 u_2}{\partial x^2}\end{pmatrix}
      + \begin{pmatrix}0 & -2\\2 & 0\end{pmatrix}
        \begin{pmatrix}\frac{\partial^2 u_1}{\partial x\partial y}\\\frac{\partial^2 u_2}{\partial x \partial y}\end{pmatrix}
        - \begin{pmatrix}1 & 0\\0 & 1\end{pmatrix}
          \begin{pmatrix}\frac{\partial^2 u_1}{\partial y^2}\\\frac{\partial^2 u_2}{\partial y^2}\end{pmatrix} =
          \begin{pmatrix}0\\0\end{pmatrix}
  \end{equation*}
  The determinant associated with the quadratic form is
  \begin{equation*}
    \det\left(
    \begin{pmatrix}1 & 0\\0 & 1\end{pmatrix} \lambda^2
      + \begin{pmatrix}0 & -2\\2 & 0\end{pmatrix} \lambda
        - \begin{pmatrix}1 & 0\\0 & 1\end{pmatrix}
          \right) =
          \det \begin{pmatrix}\lambda^2 -1 & -2\lambda\\2\lambda & \lambda^2-1\end{pmatrix} \neq 0
  \end{equation*}
  for any $\lambda\in\real$.
  So, system (\ref{eq:4}) is elliptic in the sense of \cite{2,5}.
  
}
\begin{equation}
  \left\{
    \begin{aligned}
      \frac{\partial^2 u_1}{\partial x^2} - 
      2 \frac{\partial^2 u_2}{\partial x \partial y} -
      \frac{\partial^2 u_1}{\partial y^2} & = 0, \\
      \frac{\partial^2 u_2}{\partial x^2} + 
      2 \frac{\partial^2 u_1}{\partial x \partial y} -
      \frac{\partial^2 u_2}{\partial y^2} & = 0,
    \end{aligned} \right.
  \label{eq:4} 
\end{equation}
when domain $D^-$ is the exterior of unit circle $|z| \le 1$.

This system can be written in the following form:
\begin{equation}\label{eq:26}
  \frac{\partial^2 w}{\partial \bar z^2} = 0
\end{equation}
where 
\begin{gather*}
  2\dfrac{\partial \omega}{\partial \bar z} = \dfrac{\partial \omega}{\partial
    x} + i \dfrac{\partial \omega}{\partial y}\\
  w = u_1(x,y) + i u_2(x,y)
\end{gather*}

From (\ref{eq:26}) one can conclude immediately that the regular
solution to system (\ref{eq:4}) can be given in the following way:
\begin{equation*}
  w = u_1(z) + i u_2(z) = u_1(x,y) + i\, u_2(x,y) = \bar z \varphi(z) + \psi(z)
\end{equation*}
for any arbitrary analytic function $\varphi(z)$ and $\psi(z)$ of
variable $z$ in the region $D^-$ \cite{1}.

\paragraph{Homogeneous Dirichlet problem.}
It is easy to check that homogeneous Dirichlet problem
\begin{equation*}
  \begin{gathered}
  u_1(t) = 
  u_2(t) = 0,
  \quad t \in S,  \qquad
  P = Q = \begin{pmatrix}
    0 & 0 \\
    0 & 0
  \end{pmatrix}, \quad
  R = \begin{pmatrix}
    1 & 0\\
    0 & 1
  \end{pmatrix}\\
  \det(P + i Q) = 0, \qquad \det(R) = 1,  \qquad f(t) = (f_1(t),f_2(t)), \qquad x + i y = t = e^{i\theta}
  \end{gathered}
\end{equation*}
for system (\ref{eq:4}), in region $D^-$, has an infinite set of
linearly independent regular solutions,
\begin{equation*}
  \omega_k(z) = u_{1k}(z) + i u_{2k}(z) = \bar z z^{-k} - z^{-(k+1)},\qquad
  k \ge 1,\quad, z\in D^-
\end{equation*}
bounded at infinity.
    

\paragraph{Homogeneous Neumann problem.}

On the other hand, let us consider the following homogeneous Neumann
problem in the above region $D^-$
\begin{equation*}
  \begin{gathered}
  \left.\begin{aligned}
    \frac{\partial u_{1}}{\partial \nu} &= 0\\
    \frac{\partial u_{2}}{\partial \nu} &= 0
  \end{aligned}\right\}\qquad  
  P = \begin{pmatrix}
    \cos \widehat{\nu x} & 0 \\
    0 & \cos \widehat{\nu x}
  \end{pmatrix}, \qquad
  Q = \begin{pmatrix}
    \cos \widehat{\nu y} & 0 \\
    0 & \cos \widehat{\nu y}
  \end{pmatrix}, \\
  R = \begin{pmatrix}
    0 & 0 \\
    0 & 0
  \end{pmatrix},  \qquad
  f(t)=(f_1(t),f_2(t)) \qquad
  t = e^{i\theta}
  \end{gathered}
\end{equation*}
In this case, $\cos \nu y = \sin \nu x$
\begin{multline*}
  P + iQ
  = \begin{pmatrix}\cos\widehat{\nu x} & 0 \\0 & \cos\widehat{\nu x}\end{pmatrix} + i \begin{pmatrix}\sin\widehat{\nu x} & 0 \\ 0 & \sin\widehat{\nu x}\end{pmatrix}\\
      = \begin{pmatrix}\cos\widehat{\nu x} + i  \sin\widehat{\nu x}& 0 \\0 & \cos\widehat{\nu x} + i \sin\widehat{\nu x}\end{pmatrix}
\end{multline*}
\begin{multline*}
        \det(P+iQ)
        = (\cos\widehat{\nu x} + i  \sin\widehat{\nu x})^2 
        = \cos^2 \widehat{\nu x} -  \sin^2\widehat{\nu x} + 2 i \cos \widehat{\nu x} \sin\widehat{\nu x}\\
        = \cos 2\widehat{\nu x} + i  \sin2\widehat{\nu x} = \exp(i 2\widehat{\nu x}) \neq 0
\end{multline*}
where $\theta = 2\widehat{\nu x}$.



This problem, as linearly independent regular solutions, has functions
\begin{equation*}
  \begin{aligned}
    \omega_0 &= 1,\\
    \omega_k(z) &= u_{1k}(z) + i\, u_{2k}(z) = \bar z \frac 1{z^k} - \frac{k-1}{k+1} \frac 1{z^{k+1}}, 
    \qquad k=1,\ldots.
  \end{aligned}
\end{equation*}



The above example shows that the fulfillment of condition (\ref{eq:3})
does not always guarantee the normal solvability of problem
(\ref{eq:4}), (\ref{eq:2}) in the sense of Noether.

\paragraph{Neumann problem.}

One should not assume, however, that for system (\ref{eq:4}) there is
no problem of type (\ref{eq:2}) with the property of normal
solvability. To this end, let us assume that in boundary-value
condition (\ref{eq:2}) matrices $P$, $Q$, $R$ are selected as follows
\begin{equation*}
  \begin{gathered}
    P =  \begin{pmatrix} 0 & 0 \\ 1 &  0  \end{pmatrix}, \qquad
    Q = \begin{pmatrix}  0 & 0 \\ 0 &  -1 \end{pmatrix},\qquad
    R =  \begin{pmatrix} 1 & 0 \\ 0 &  0  \end{pmatrix},\\
    \det(P + i Q) = \det\begin{pmatrix}0&0\\1&-i\end{pmatrix} =0,\qquad
    f(t) = (f_1(t),f_2(t))\qquad
    t = e^{i\theta}
  \end{gathered}
\end{equation*}
In the considered case the search regular solution will be as follows,
\begin{multline*}
  \omega(z) = - \frac{\bar z}{4 \pi i} \int_0^{2 \pi}
  \frac{t+z}{t(t-z)} f_2(t) \, dt + \frac 1{4 \pi i} \int_0^{2\pi}
  \frac{t+z}{zt(t-z)} f_2(t)\, dt\\ - \frac 1{\pi i} \int_0^{2\pi}
  \frac {f_1(t)}{t-z} \, dt + \frac 1{2 \pi i} \int_0^{2\pi}
  \frac{f_1(t)}t\, dt + i K \qquad z \in D^-
\end{multline*}
\begin{equation}
  w(z) = u_1(z)+ u_2(z) =   \bar z \varphi(z) + \psi(z),
  \qquad z\in D^-, t\in S\label{eq:27}
\end{equation}
This solution fulfills (\ref{eq:4}) and boundary-value condition
(\ref{eq:2}) when $\int_0^{2\pi} f_2(\theta)\, d\theta = 0$,
$t=e^{i\theta}$.




Unlike the Neumann problem considered above, condition (\ref{eq:3})
might be fulfilled and problem (\ref{eq:4}), (\ref{eq:2}) is normally
solvable. To confirm this assumption we can consider problem
(\ref{eq:4}), (\ref{eq:2}) in the assumption that
\begin{equation}
  \begin{gathered}
    P = \begin{pmatrix} 1&0\\0&1\end{pmatrix},\quad 
      Q = \begin{pmatrix} 0&-1\\-1&0\end{pmatrix},\quad
        R = \begin{pmatrix} 0&0\\0&0\end{pmatrix}, \\
          \det(P+iQ) = \det\begin{pmatrix}1&-i\\-i &1\end{pmatrix} = 2,\qquad
          f(t)= (f_1(t),f_2(t)), \qquad t=e^{i\theta}
  \end{gathered}
          \label{eq:6}
\end{equation}
On the basis of (\ref{eq:27}) and (\ref{eq:6}), let us write the
boundary-value condition (\ref{eq:2}) in the form
\begin{equation}
  2 \Re \frac{\partial \omega}{\partial \bar t} = 
  2 \Re \varphi(t) = f_1(t), \qquad t\in S, \label{eq:7}
\end{equation}
\begin{equation}2 \Im \frac{\partial \omega}{\partial t} = 
  2 \Im\left[\bar t \varphi'(t) + \psi'(t)\right] = f_2(t),
  \quad t\in S. \label{eq:8}
\end{equation}


Taking into the account (\ref{eq:27}) one can write
\begin{equation}
  \varphi(z) = -\frac{1}{4\pi i} \int_0^{2\pi} \frac{t+z}{t(t-z)} f_2(t)\, dt
  \label{eq:29}
\end{equation}
As the point $z=x+iy$ from the region $D^-$ approaches the boundary
point $t\in S$, on the basis of considering the class of analytical
functions in the domain $D^-$ that vanish at infinity, from
(\ref{eq:7}) and by the Morera and Cauchy theorems, one has
\begin{equation*}
  \varphi(t) = 0 = \int_0^{2\pi} f_1(t)\, dt \Leftrightarrow f_1(t)
  \text{ is analytical.}
\end{equation*}
Therefore, in the class of analytical functions in the domain $D^-$
vanishing at infinity, problem (\ref{eq:7}) is solvable if and only if
\begin{equation}
  \int_0^{2\pi} f_1(t)\, dt = 0, \qquad t=e^{i\theta}
  \label{eq:9}
\end{equation}
Its solution is unique, being given by
\begin{equation*}
  \varphi(z) = -\frac{1}{4\pi i} \int_0^{2\pi} \frac{t+z}{t(t-z)} f_1(t)\, dt
\end{equation*}

On the other hand, in the case of (\ref{eq:8}) we obtain
\begin{equation}
  \int_0^{2\pi}f_2(\theta)\, d\theta = 0, \label{eq:11}
\end{equation}
and when condition (\ref{eq:11}) is fulfilled, the analytical function
$\psi(z)$ in $D^-$, vanishing at infinity, can be defined by formula
\begin{equation}
  \psi'(z) = - \frac{\varphi'(z)}z - \frac 1{4\pi} 
  \int_0^{2\pi}\frac{t+z}{t(t-z)}f_2(t)\,dt.
  \label{eq:12}
\end{equation}

Function $\psi(z)$ can be determined from (\ref{eq:12}) as a result of
integration up to the arbitrary complex constant. Therefore, it is
possible to conclude that for the solvability of non-homogeneous
problem (\ref{eq:4}), (\ref{eq:2}), when matrices $P$, $Q$, $R$ are
given by formulas (\ref{eq:6}), it is necessary and sufficient that
functions $f_1$ and $f_2$ in the right parts of the boundary-value
condition (\ref{eq:2}) satisfy equalities (\ref{eq:9}) and
(\ref{eq:11}). When these equalities are fulfilled the solution of the
problem can be defined up to the arbitrary additive constant.


\paragraph{Conclusion.}

The above examples show the importance of the definition of the
classes of elliptic systems for which the Poincaré problem, in the
case of exterior domains, is normally solvable according to Noether.


\section{Decomposable system}

Consider the linear system of equations in partial derivatives of type
(\ref{eq:1}) in the case when the characteristic
determinant
\begin{equation*}
  \det(A \lambda^2 + 2 B \lambda + C) \neq 0 \label{eq:13}
\end{equation*}
for any $\lambda\in\real$. It is known that such systems are called
elliptic \cite{2}. For reference, the sixth chapter of \cite{25} is
devoted entirely to de Poincare problem for second order linear
elliptic systems.

\begin{defn}\label{def:3}
  System (\ref{eq:1}) is called decomposable \cite{2}, if $A$, $B$,
  $C$ are diagonal matrices, i.e. $A = (a_{jk})$, $B = (b_{jk})$, $C =
  (c_{jk})$,
  \begin{equation*}
    a_{jk} = 
    \begin{cases}
      a_j, & k = j, \\
      0, & k\neq j,
    \end{cases}\quad
    b_{jk} = 
    \begin{cases}
      b_j, & k = j, \\
      0, & k\neq j,
    \end{cases}\quad
    c_{jk} = 
    \begin{cases}
      c_j, & k = j, \\
      0, & k\neq j,
    \end{cases}
    \quad j,k = 1,\ldots,n.
  \end{equation*}
\end{defn}
  
Let $A$, $B$, $C$ be constant diagonal matrices, meaning that the
considered system is decomposed into $n$ equations.
\begin{equation}
  a_j \frac{\partial^2 u_j}{\partial x^2} + 
  2 b_j \frac{\partial^2 u_j}{\partial x \partial y} +
  c_j \frac{\partial^2 u_j}{\partial y^2}  = 0, \qquad j = 1, \ldots, n. \label{eq:14}
\end{equation}

In these assumptions the ellipticity of system (\ref{eq:1}) means that
all quadratic forms
\begin{equation*}
  Q_j(\lambda_1, \lambda_2) = 
  a_j \lambda_1^2 + 2 b_j \lambda_1 \lambda_2 + c_j \lambda_2^2,
  \qquad j = 1,\ldots,n,
\end{equation*}
are positively definite.

\begin{defn}\label{def:4}
  The regular solution of equation (\ref{eq:14}) at each
  $j=1,\ldots,n$ in $D^-$ will be called the bounded in $D^-$
  twice-continuously differentiable function $u_j(x,y)$, satisfying
  this equation in each finite point in $D^-$.
\end{defn}

Likewise, it is possible to verify directly that function $u_j(x,y)$
defined by formula
\begin{equation}
  u_j(x,y) = \Re \varphi_j(z_j) \label{eq:15}
\end{equation}
where $\varphi_j(z_j)$ is the arbitrary bounded analytical function in
$D^-$ domain of complex variable
\begin{equation*}
  \begin{aligned}
    z_j &= \frac 1{\sqrt{a_j}} x + 
    i \left(\frac{\sqrt{a_j}}{\delta_j} y - \frac{b_j}{\delta_j\sqrt{a_j}}x\right),\\
    \delta_j^2 &= a_j c_j - b_j^2, 
  \end{aligned} \label{eq:16}
\end{equation*}
gives regular solutions of equation (\ref{eq:14}) in $D^-$. The proof
of this proposition is available in the appendix of this paper.

It is evident that solution $u_j(x,y)$ of equation (\ref{eq:14}),
regular in $D^-$, function $\varphi_j(z_j)$ can be defined up to an
arbitrary imaginary constant. Without loss of generality, it can be
assumed that function $\varphi_j(z_j)$ in formula (\ref{eq:15})
satisfies condition
\begin{equation*}
  \Im \varphi_j(\infty) = 0.
\end{equation*}

The following theorem will be taking into account below:


\begin{thm}[Representation theorem, Vekua \cite{3}]\label{thm:1}
  For each bounded analytical function $\varphi_j(z_j)$ of
  $C^{1,h}(D^- \cup S)$ above considered\footnote{The integral in
  formula (\ref{eq:18}) represents the integral at one circumvention
  of point $t$ around $S$ from the fixed point $z_0=0$ on $S$ in the
  positive direction
  

  } there exists unique real function
  $\mu_j(s)$ of the class $C^{0,h} (S)$ so that
  \begin{equation}
    \varphi_j(z_j) = \int_S\ln\left(1-\frac{z_j}{\tau}\right) \mu_j(\tau)\, 
    d s_{\tau}. \label{eq:18}
  \end{equation}
\end{thm}

Below, the Poincaré problem for decomposable system (\ref{eq:14}) is
investigated and formulated as follows:
\begin{thm}[Poincare problem for decomposable systems]\label{thm:2}
  Let's consider the decomposable system with boundary-value
  conditions
  \begin{equation}
    \begin{gathered}
      a_j \frac{\partial^2 u_j}{\partial x^2} + 
      2 b_j \frac{\partial^2 u_j}{\partial x \partial y} +
      c_j \frac{\partial^2 u_j}{\partial y^2}  = 0, \qquad j = 1, \ldots, n.\\
      P \frac{\partial u}{\partial x} + Q \frac{\partial
        u}{\partial y} = f(t),
      \qquad x + i y = t \in S
    \end{gathered}
    \label{eq:17}
  \end{equation}
  where $a_j, b_j, c_j$ are constant values, $j=1,\ldots,n$, and $P,
  Q, R$ given real $n\times$ matrices on $S$ satisfying H\"older
  conditions, $f(t)=(f_1(t),\ldots,f_n(t))$ given real vector on $S$
  satisfying H\"older's condition.

  The solution to this problem is equivalent to finding the solution
  to a system of singular integral equation.
  
\end{thm}

\begin{proof}
  The proof is based on Theorem \ref{thm:1}.

  Returning to formula (\ref{eq:15}), we can write boundary-value
  condition (\ref{eq:17}) as
  \begin{multline}
    \Re\sum_{j=1}^n
    \left[
      p_{jk}\left(
      \frac 1{\sqrt{a_j}} - i \frac{b_j}{\delta_j\sqrt{a_j}}
      \right) + i q_{jk} \frac{\sqrt{a_j}}{\delta_j} 
      \right] \varphi_j'(t) = f_k(t), \\
    k = 1,\ldots,n, \quad t(s)\in S. 
    \label{eq:19}
  \end{multline}
  
  Taking symbols $\alpha = (\alpha_{jk})$, $\beta=(\beta_{jk})$,
  \begin{align}
    \alpha_{jk}(t) &= \Re\left[ p_{jk}(t) 
      \left(
      \frac 1{\sqrt{a_j}} - i \frac{b_j}{\delta_j\sqrt{a_j}}
      \right) + i q_{jk}(t) \frac{\sqrt{a_j}}{\delta_j} \right] 
    \pi i \bar t', \label{eq:20}\\
    \beta_{jk}(t) &= \Im\left[ p_{jk}(t) 
      \left(
      \frac 1{\sqrt{a_j}} - i \frac{b_j}{\delta_j\sqrt{a_j}}
      \right) + i q_{jk}(t) \frac{\sqrt{a_j}}{\delta_j} \right]
    i \bar \tau', \label{eq:21}
  \end{align}
  and based on (\ref{eq:15}), (\ref{eq:19}) and (\ref{eq:18}), we
  obtain
  \begin{equation}
    T_\mu \equiv \alpha(t) \mu(t) - \beta(t) \int_S\frac{\mu(\tau)}{\tau-t}\,d\tau +
    \int_SK(t,\tau)\mu(\tau)\,d\tau = f(t), \label{eq:22}
  \end{equation}
  where $\mu=(\mu_1,\ldots,\mu_n)$ is search vector
  \begin{equation}
    K(t,\tau) =
    \frac{K^1(t,\tau)-K^1(t,t)}{\tau-t}, \label{eq:23}
  \end{equation}
  and $n\times n$ matrix $K(t,\tau)$ is given by
  \begin{multline}
    K^1_{jk}(t,\tau) = \bar\tau'(\tau-t) 
    \Re\left\{\left[p_{jk}(t)\left(
      \frac 1{\sqrt{a_j}} - i \frac{b_j}{\delta_j\sqrt{a_j}}
      \right) + iq_{jk}(t) \right]
    \left[\frac 1{\tau -t} - K_{0j}(t,\tau)\right]\right\},\\
    K_{0j}(t,\tau) \equiv  -\frac{\bar\tau_j'}{t_j} - 
    \frac{d\bar\tau_j}{d s_{\tau_j}} \frac 1{\tau_j-t_j}  -
    \frac{d\,\bar\tau}{d s_\tau} \frac 1{\tau-t}, \qquad j,k = 1,\ldots,n
    \label{eq:24}
  \end{multline}
  
  Due to (\ref{eq:23}) and (\ref{eq:24}) the last term in the left
  part of equality (\ref{eq:22}) represents a compact matrix integral
  operator and the second term a matrix singular integral operator, in
  which the integral is understood in the sense of Cauchy's principal
  value. Likewise, due to (\ref{eq:20}) and (\ref{eq:21}) matrix
  coefficients $\alpha(t)$ and $\beta(t)$ satisfy H\"older's
  condition.

  Therefore, (\ref{eq:22}) is a system of singular integral equations
  equivalent to problem (\ref{eq:17}).

  It is known that system (\ref{eq:22}) is normally solvable if the
  conditions
  \begin{equation}
    \left\{
    \begin{aligned}
      \det[\alpha(t) - \pi i \beta(t)] &\neq 0\\
      \det[\alpha(t) + \pi i \beta(t)] &\neq 0
    \end{aligned}
    \right.\label{eq:25}
  \end{equation}
  are fulfilled everywhere on $S$.

  The index of this system is calculated by formula \cite{2,9}:
  \begin{equation*}
    \kappa = \frac 1{2 \pi} 
    \left[
      \arg\frac {\det(\alpha + i \pi \beta)}{\det(\alpha - i \pi \beta)} 
      \right]_S
  \end{equation*}
  where $[]_S$ means the increment of the function in brackets at one
  circumvention of point $t$ around of curve $S$ from the fixed point
  $z_0=0$ on $S$ in the positive direction.


  
  When fulfilling conditions (\ref{eq:25}) the well-known Noether
  theorems are used for the system of singular integral equations
  (\ref{eq:22}):
  \begin{enumerate}
  \item the corresponding (\ref{eq:22}) homogeneous system $T\mu_0=0$
    and adjoin homogeneous system $T^*\mu_0=0$ have no more than
    finite number $l$ and $l'$ of linearly independent solutions,
    \begin{equation*}
      l-l'=\kappa,
    \end{equation*} 
  \item for the solvability of non-homogeneous system (\ref{eq:22}) it
    is necessary and sufficient that conditions
    \begin{equation*}
      \int_Sf(t)\mu_*^{(k)}(t)\, dt = 0,\qquad k=1,\ldots,l'.
    \end{equation*}
    are fulfilled, where vectors $\{\mu_*^{(k)}\}$, $k=1,\ldots,l'$,
    represent all linearly independent solutions of ad-joint
    homogeneous system $T^*\mu_*=0$.
  \end{enumerate}
\end{proof}
   
On the basis of equivalence between problem (\ref{eq:17}) and the
system of singular integral equations (\ref{eq:22}), it is possible to
conclude that for the solvability of this problem for arbitrary $f\in
C^{0,h}(S)$ it is necessary and sufficient that $l'=0$, and the
general solution to equation (\ref{eq:22}) has the following form
\begin{equation*}
  \mu = \sum_{k=1}^l \beta_k \mu_k + \mu_0
\end{equation*}
where $\mu_k$ represents all linearly independent solutions of
equation $T_{\mu_0}$. $\beta_k$ are real arbitrary constants and
$\mu_0$ partial solution of the same equation.

If $\kappa=0$ and the homogeneous problem corresponding to
(\ref{eq:17}) has only a trivial solution, then the non-homogeneous
problem (\ref{eq:17}) has a solution which is the only one.

\section{Conclusions}

In this work, we have studied the exterior boundary value Poincaré
problem for an elliptic system posed in equation (\ref{eq:1}) that
satisfy the boundary conditions given by (\ref{eq:2}). We have offered
some examples to show that the requirement of uniform ellipticity of
system (\ref{eq:1}) does not guarantee the normal solvability
(according to Noether). Nevertheless, we prove that, under certain
conditions, it is possible to select from system (\ref{eq:1}) some
decomposable types of partial differential equations (given by
equation (\ref{eq:14})) that are effectively solvable in the sense of
Noether. For those systems, the Poincaré problem (equation
(\ref{eq:17})) is reduced by means of Theorem \ref{thm:1} to an
equivalent system of singular integral equations.

The importance of integral equation methods in the solution of certain
types of boundary value problems is universally accepted. From a
practical point of view, the benefit may not always be very relevant
when interior problems are concerned, but for exterior problems, where
the region of interest is the infinite extent, an integral equation
formulation may be virtually indispensable and our results are able to
provide added value on this topic.

\section*{Appendix}

The solution presented in equation \ref{eq:15} can be reformulated
through a variable substitution:
\begin{gather*}
  x_j(x,y) = \frac 1{\sqrt{a_j}} x, \qquad
  y_j(x,y) = \frac{\sqrt{a_j}}{\delta_j} y - \frac{b_j}{\delta_j\sqrt{a_j}} x, \qquad
  \delta_j^2 = a_j c_j - b_j^2, \\
  u_j(x,y) = \Re \varphi_j(z_j) = w_j(x_j,y_j), \quad z_j=x_j+iy_j
\end{gather*}
we canexpress it as:
\begin{equation*}
  u_j(x,y) = w_j(x_j(x,y), y_j(x,y))
\end{equation*}
where $w_j$ is a harmonic function corresponding to the real part of a
analytical complex function $\varphi_j$.

By applying the chain rule, we can compute the following derivatives:
\begin{align*}
  \frac{\partial^2 u_j}{\partial x^2}
  &= \frac{\partial^2 w_j}{\partial x_j^2}\frac{1}{a_j}
  - 2\frac{\partial^2 w_j}{\partial x_j\partial y_j}\frac{b_j}{\delta_ja_j}
  + \frac{\partial^2 w_j}{\partial y_j^2}\frac{b_j^2}{\delta_j^2a_j}
  \\
  \frac{\partial^2 u_j}{\partial y^2}
  &= \frac{\partial^2 w_j}{\partial y_j^2}\frac{a_j}{\delta_j^2}
  \\
  \frac{\partial u_j}{\partial x\partial y}
  &= \frac{\partial^2 w_j}{\partial y_j^2} \frac{-b_j}{\delta_j^2}
  + \frac{\partial^2 w_j}{\partial x_j\partial y_j} \frac{1}{\delta_j}
\end{align*}
By substituing the previous expressions into equation (\ref{eq:14}):
\begin{equation*}
  a_j \frac{\partial^2 u_j}{\partial x^2} + 
  2 b_j \frac{\partial^2 u_j}{\partial x \partial y} +
  c_j \frac{\partial^2 u_j}{\partial y^2}  = 0,
\end{equation*}
after some straightforward algebraic operations, it is equivalent to
\begin{equation*}
  \frac{\partial^2 w_j}{\partial x_j^2}
  + \frac{\partial^2 w_j}{\partial y_j^2}
  = 0
\end{equation*}
that holds because $w_j$ is a harmonic function.



%

\section*{Declarations}
\paragraph{Ethical approval}\ \\
Not applicable.

\paragraph{Availability of data and materials}\ \\
Data sharing is not applicable to this article as no datasets were
generated or analyzed during the current study.

\paragraph{Competing interests}\ \\
The authors declare that they have no competing interests.

\paragraph{Funding}\ \\
The Authors received no specific funding for this work.

\paragraph{Authors' contribution}\ \\
All the authors contributed equally to all sections of the manuscript.

\paragraph{Acknowledgments}\ \\
Not applicable.


\end{document}